\documentclass{article}

\usepackage{amssymb}

\newcommand{\proof}{{\bf Proof:  }}
\newcommand{\remark}{{\bf Remark:  }}
\newcommand{\remarks}{{\bf Remarks:  }}
\newcommand{\example}{{\bf Example:  }}
\newcommand{\examples}{{\bf Examples:  }}
\newcommand{\dimv}{\underline{\dim}}

\newcommand{\hb}{\newline\hspace*{\fill}$\blacksquare$}

\newcommand{\od}{\overline{d}}
\newcommand{\on}{\overline{n}}
\newcommand{\oi}{\overline{I}}
\newcommand{\oq}{\overline{Q}}

\newcommand{\Imm}{{\rm Im}}
\newcommand{\Ker}{{\rm Ker}}

\newcommand{\ses}[3]
{\mbox{$0 \rightarrow #1 \rightarrow #2 \rightarrow #3 \rightarrow 0$}}

\newtheorem{theorem}{Theorem}[section]
\newtheorem{lemma}[theorem]{Lemma}
\newtheorem{definition}[theorem]{Definition}
\newtheorem{proposition}[theorem]{Proposition}
\newtheorem{corollary}[theorem]{Corollary}

\begin{document}

\parindent0pt

\title{\bf Framed quiver moduli, cohomology, and quantum groups}

\author{Markus Reineke\\ Mathematisches Institut\\ Universit\"at M\"unster\\ D - 48149 M\"unster, Germany\\ e-mail: reinekem@math.uni-muenster.de}

\date{}

\maketitle

\begin{abstract}

Framed quiver moduli parametrize stable pairs consisting of a quiver representation and a map to a fixed graded vector space. Geometric properties and explicit realizations of framed quiver moduli for quivers without oriented cycles are derived, with emphasis on their cohomology. Their use for quantum group constructions is discussed.
\end{abstract}

\section{Introduction}

In formulating interesting moduli problems for objects in an abelian category, one is usually confronted with the problem of either restricting the class of objects to be parametrized to certain rigid, or stable, objects, or to endow them with some additional structure, which rigidifies them; see \cite{M} for a thorough discussion of this problem.\\[1ex]
Moduli spaces for representations of quivers, introduced in \cite{Ki}, provide an interesting testing ground for techniques of moduli theory. This class of moduli spaces relates to many moduli problems arising in Algebraic Geometry (in particular, in connection to vector bundles). It is also interesting in itself, both in view of the classification problem of quiver representations, and, via the Hall algebra approach of \cite{Ri}, to quantum group realizations (see, for example, \cite{Re}).\\[1ex]
While the first approach to the problem of formulating interesting moduli problems for representations of quivers was used in \cite{Ki} by introducing a concept of (semi-)stable representations, the second one was first considered in \cite{Na}, as a preliminary step in the construction of Nakajima quiver varieties. In the latter case, the resulting moduli spaces are called framed quiver moduli; their construction will be recalled in section \ref{framedmoduli}.\\[2ex]
The aim of this note is to study the geometry of framed quiver moduli in the case of quivers without oriented cycles.
This study is based on two realizations, the first in terms of the ordinary quiver moduli of \cite{Ki}, the second in terms of Grassmannians of subrepresentations; they are developed in section \ref{s3}.\\[1ex]
In section \ref{geometry}, after first characterising non-vanishing of framed quiver moduli (Proposition \ref{critnon0}), we construct an embedding into a product of Grassmannians (Proposition \ref{propemb}). This gives rise to an explicit description of framed quiver moduli as iterated Grassmann bundles (Theorem \ref{chain}).\\[1ex]
Most of the results of the present paper arose out of the attempt to study the cohomology of these framed quiver moduli, with the aim of giving more precise versions of some general statements (\cite{ES,KW,Re}). The paper's main result in this direction, Theorem \ref{coho}, in fact gives very detailed infomation even on the Chow ring of framed quiver moduli, in terms of Schubert calculus (see \cite{Fu}). This follows from the realization of framed moduli as a tower of Grassmann bundles.\\[1ex]
In section \ref{orb}, a natural action of a linear algebraic group on framed quiver moduli is defined. The main result of this section, Theorem \ref{bij}, gives a description of the orbit structure of this action in terms of the orbit structure of representation varieties of quivers with respect to the base change group, which is rather well-understood.\\[1ex]
Finally, section \ref{qgc} relates the subject of framed quiver moduli to quantum group constructions. The point of view advocated there is that framed quiver moduli provide a candidate for a generalization of the apporach of \cite{BLM} to quantum groups of arbitrary (symmetric) Kac-Moody type. The main result, Theorem \ref{corsurj}, at least provides a realization of quotients of modified quantized Borel algebras in terms of a convolution algebra defined using framed quiver moduli.\\[2ex]
All the results of this note strongly depend on the hypothesis that the underlying quiver has no oriented cycles, since this allows for inductive constructions. For quivers with oriented cycles, the geometry of framed moduli is expected to be of a quite different nature (for the case of the one-loop quiver, resp.~the multiple-loop quiver, see \cite{Ge,LBR} and \cite{noncommhilb}, respectively).\\[3ex]
{\bf Acknowledgments:} I would like to thank K.~Bongartz for allowing me to use the results of \cite{Bo}, A.~Buch for discussions on the Chern class formula \cite{BF}, and S.~K\"onig for discussions on geometric realizations of $q$-Schur algebras.

\section{Definition of framed quiver moduli}\label{framedmoduli}

Let $Q$ be a finite quiver with set of vertices $I$. The arrows in $Q$ will be written as
$\alpha:i\rightarrow j$ for $i,j\in I$. Throughout the paper we assume that $Q$ does not contain
oriented cycles. Therefore, we can consider the partial order induced on $I$ by the arrows, i.e. defined by $i<j$ if there exists an arrow $\alpha:i\rightarrow j$.\\[1ex]
We denote by ${\rm rep}_kQ$ the abelian $k$-linear category of finite dimensional $k$-representations of $Q$. The homomorphisms between objects $M$ and $N$ of ${\rm rep}_kQ$ will be denoted by ${\rm Hom}_Q(M,N)$. For all basic notions on representations of quivers, the reader is referred to \cite{ARS,RiB}.\\[1ex]
Let ${\bf Z}I$ be the free abelian group in $I$; its elements will be written as
$d=\sum_{i\in I}d_ii$. Let ${\bf N}I$ be the subsemigroup of ${\bf Z}I$ of nonnegative linear
combinations of the elements of $I$; this will be viewed as the space of dimension types of
representations of $Q$. We have a non-symmetric inner product (the Euler form) on ${\bf Z}I$
defined by
$$\langle i,j\rangle=\delta_{i,j}-\#\{\mbox{arrows from $i$ to $j$}\}$$
for $i,j\in I$. It has the property
$$\langle \dimv M,\dimv N\rangle=\dim{\rm Hom}_Q(M,N)-\dim{\rm Ext}^1_Q(M,N)$$
for any two representations $M$ and $N$ of $Q$. Later on, we will use the symmetrization $(d,e)=\langle d,e\rangle+\langle e,d\rangle$. We also consider the standard inner product $i\cdot j=\delta_{i,j}$ on ${\bf N}I$.\\[1ex]
Let $k$ be an algebraically closed field, and let $d,n$ be dimension types, that is, elements of ${\bf N}I$. For $i\in I$, let $M_i$ (resp.~$V_i$) be a $k$-vector
space of dimension $d_i$ (resp.~$n_i$). In the following, we will frequently consider $M=\bigoplus_{i\in I}M_i$ as an $I$-graded vector space. In particular, we consider $I$-tuples of linear maps like $(f_i:M_i\rightarrow V_i)_{i\in I}$ as linear maps $f:M\rightarrow V$ of $I$-graded vector spaces.\\[1ex]
We define the variety of $k$-representations
of $Q$ of dimension type $d$ by
$$R_d=R_d(Q)=\bigoplus_{\alpha:i\rightarrow j}{\rm Hom}_k(M_i,M_j),$$
where the sum runs over all arrows in $Q$. The group $G_d=\prod_{i\in I}{\rm GL}(M_i)$ acts on $R_d
$ via
$$g\cdot M=(g_i)_{i\in I}\cdot(M_\alpha)_{\alpha:i\rightarrow j}=
(g_jM_\alpha g_i^{-1})_{\alpha:i\rightarrow j}.$$
The variety $R_d$ is an affine space of dimension $\dim R_d=\sum_{\alpha:i\rightarrow j}d_id_j$,
and the group $G_d$ is reductive and algebraic of dimension $\dim G_d=\sum_{i\in I}d_i^2$. In particular, note that
$$\dim G_d-\dim R_d=\langle d,d\rangle.$$
By definition, the orbits ${\cal O}_M$ of $G_d$ in $R_d$ are in bijection with the isomorphism
classes $[M]$ of $k$-representations of $Q$ of dimension type $d$. An orbit ${\cal O}_M$ is closed in $R_d$ if and only if the representation $M$ is semisimple. Since $Q$ has no oriented cycles, there exists a unique such representation (up to isomorphism) for each dimension vector. Thus, $G_d$ acts on $R_d$ with a unique closed orbit, for each $d\in{\bf N}I$.\\[1ex]
Now we consider the extended representation variety
$$R_{d,n}=R_d\times\bigoplus_{i\in I}{\rm Hom}_k(M_i,V_i).$$
Again, it carries a natural action of $G_d$ given by
$$g\cdot(M,(f_i)_{i\in I})=(g\cdot M,(f_ig_i^{-1})_{i\in I}),$$
extending the action of $G_d$ on $R_d$.\\[1ex]
We introduce a notion of stability for points of $R_{d,n}$:

\begin{definition}\label{ds} A point $(M,f=(f_i)_{i\in I})$ of $R_{d,n}$ is called stable if there is no 
non-zero subrepresentation
$U$ of $M$ which is contained in ${\rm Ker}f\subset M$. The set of all stable points of $R_{d,n}$ is denoted by $R_{d,n}^s$.
\end{definition}

Now we can define the framed quiver moduli:

\begin{theorem}[\cite{Na}]\label{deffm} There exists a geometric quotient ${\cal M}_{d,n}={\cal M}_{d,n}(Q)$ of $R_{d,n}^s$ by $G_d$.
\end{theorem}

The spaces ${\cal M}_{d,n}$ are called framed quiver moduli. They can be defined explicitly in terms of semi-invariants on the extended representation varieties as follows (see \cite{Na}): define a character $\chi$ on the group $G_d$ by
$$\chi((g_i)_{i\in I}):=\prod_{i\in I}\det(g_i).$$
Let $A=k[R_{d,n}]$ be the coordinate ring of the affine variety $R_{d,n}$. For $n\in{\bf N}$, let $A^{G_d,\chi^n}$ be the space of semi-invariant functions with respect to the character $\chi^n$, that is, regular functions $f\in A$ such that
$$f(g\cdot x)=\chi(g)^n\cdot f(x)\mbox{ for all $g\in G_d$ and $x\in R_{d,n}$}.$$
Then
$${\cal M}_{d,n}={\rm \bf Proj}(\bigoplus_{n\in{\bf N}}A^{G_d,\chi^n}).$$
In particular, since $0$ is the unique closed $G_d$-orbit in $R_d$, we have $A^{G_d}=k$. Thus ${\cal M}_{d,n}$ is a projective variety. From Theorem \ref{deffm}, we also get immediately:

\begin{lemma} If ${\cal M}_{d,n}$ is nonempty, it is a smooth projective variety of dimension
$$\dim {\cal M}_{d,n}=d\cdot n-\langle d,d\rangle.$$
\end{lemma}

\proof By definition of ${\cal M}_{d,n}$ as a geometric quotient of $R_{d,n}^s$ by $G_d$, we have:
$$\dim {\cal M}_{d,n}=\dim R_{d,n}-\dim G_d=\sum_{\alpha:i\rightarrow j}d_id_j+\sum_{i\in I}d_in_i-\sum_{i\in I}d_i^2
=d\cdot n-\langle d,d\rangle.
$$\hb

\example As an easy (but nevertheless typical) example of framed quiver moduli, we consider the case of the quiver of type $A_1$, that is, $I=\{i\}$ is a one-element set, and there are no arrows. In this case, $R_{d,n}={\rm Hom}_k(M_i,V_i)$ for $k$-vector spaces $M_i$ and $V_i$ of dimension $d_i$ and $n_i$, respectively, and $R_{d,n}$ consists precisely of the injective homomorphisms. Thus, we see immediately that ${\cal M}_{d,n}$ is empty if $d>n$, and that ${\cal M}_{d,n}$ is isomorphic to the Grassmannian ${\rm Gr}_d(V_i)$ of $d$-dimensional subspaces of $V_i$, by associating to a point $f\in R_{d,n}^s$ the subspace $\Imm(f)$ of $V_i$.\\[2ex]
\remark The construction of framed quiver moduli has a dual version by defining the extended representation variety as
$R_d\times\bigoplus_{i\in I}{\rm Hom}_k(V_i,M_i)$ and considering stable pairs $(M,f)$ such that there is no proper subrepresentation $U$ of $M$ containing ${\rm Im} f$. The interested reader will easily dualize all the results of this paper to this situation.

\section{Two realizations of framed quiver moduli}\label{s3}

In this section, we give two realizations of the framed quiver moduli ${\cal M}_{d,n}$. The first gives an interpretation in terms of the ordinary quiver moduli introduced by A.~King in \cite{Ki}, the second in terms of Grassmannians of subrepresentations of a quiver representation.

\subsection{Realization in terms of ordinary quiver moduli - deframing}\label{real1}

We first recall A.~King`s construction of (ordinary) quiver moduli \cite{Ki}, with a slightly different notion of stability used in \cite{Re}. As before, let $Q$ be a finite quiver without oriented cycles. Choose a linear function $\Theta:{\bf Z}I\rightarrow {\bf Z}$, called a stability for $Q$, and denote by $\dim$ the linear function on ${\bf Z}I$ defined by $\dim i=1$ for all $i\in I$. Define the slope of a non-zero dimension type $d\in{\bf N}I\setminus 0$ by $\mu(d)=\frac{\Theta(d)}{\dim d}$.
Finally, define the slope $\mu(X)=\mu(\dimv X)$ of a representation $X$ of $Q$ as the slope of its dimension type.

\begin{definition} A representation $X$ of $Q$ is called semistable (resp.~stable) if $\mu(U)\leq \mu(X) \mbox{(resp. } \mu(U)<\mu(X)\mbox{)}$ for all non-zero proper subrepresentations $U$ of $X$.
Denote by $R_d^{ss}$ (resp.~$R_d^s$) the set of points of $R_d$ corresponding to semistable (resp.~stable) representations.
\end{definition} 

\begin{theorem}[\cite{Ki}]
The semistable locus $R_d^{ss}$ is an open subset of $R_d$, and the stable locus $R_d^s$ is an open subset of $R_d^{ss}$.
There exists an algebraic quotient ${\cal M}_d^{ss}$ of $R_d^{ss}$ by $G_d$. The variety ${\cal M}_d^{ss}$ is projective, of dimension $1-\langle d,d\rangle$ if it is non-empty. There exists a geometric quotient ${\cal M}_d^s$ of $R_d^s$ by $G_d$. The variety ${\cal M}_d^s$ is smooth and embeds as an open subset of ${\cal M}_d^{ss}$.
\end{theorem}

In particular, consider the case where $R_d^{ss}=R_d^s$, that is, all semistables are already stable. Then ${\cal M}_d^{ss}={\cal M}_d^s$ is a smooth projective variety, parametrizing the isomorphism classes of (semi-)stable representations of $Q$ of dimension type $d$ (since ${\cal M}_d^s$ is a geometric quotient).\\[2ex]
We introduce a new quiver $\widetilde{Q}$ with set of vertices 
$\widetilde{I}=I\dot\cup\{\infty\}$, the arrows in $\widetilde{Q}$ being those of $Q$, together with $n_i$ arrows from $i$ to $\infty$, for each $i\in I$. We also extend the dimension type to $\widetilde{d}=d+\infty\in{\bf N}\widetilde{I}$. We consider the stability $\Theta=-\infty^*$.\\[1ex]
By definition, we have a $G_d$-equivariant isomorphism
$$R_{\widetilde{d}}(\widetilde{Q})=\bigoplus_{\alpha:i\rightarrow j}{\rm Hom}_k(k^{d_i},k^{d_j})\oplus\bigoplus_{i\in I}{\rm Hom}_k(k^{d_i},k)^{n_i}\simeq R_{d,n}(Q)$$
via the identification ${\rm Hom}_k(k^{d_i},k)^{n_i}\simeq{\rm Hom}_k(k^{d_i},k^{n_i})$.

\begin{proposition} The moduli space ${\cal M}_{\widetilde{d}}^{ss}(\widetilde{Q})$, defined by the
above stability, is isomorphic to ${\cal M}_{d,n}$.
\end{proposition}

\proof Using the above identification $R_{\widetilde{d}}(\widetilde{Q})$ with $R_{d,n}(Q)$, we will show that the set of (semi-)stable points $R_{\widetilde{d}}^{ss}(\widetilde{Q})$ identifies with the set of stable points $R_{d,n}^s(Q)$ as in Definition \ref{ds}. Let $X$ be a representation of $\widetilde{Q}$ of dimension type $\widetilde{d}$, viewed as a pair $(M,f)$ consisting of a representation $M$ of $Q$ of dimension type $d$, together with an $I$-graded linear map
$f=(f_i)_{i\in I}:M\rightarrow V$. By definition, $X$ is semistable if and only if $\mu(Y)\leq\mu(X)$ for all non-zero proper subrepresentations $Y$ of $X$.\\[1ex]
We have $\mu(X)=-\frac{1}{\dim X+1}$.
If $(\dimv Y)_\infty=1$, then $\mu(Y)=-\frac{1}{\dim Y}<\mu(X)$. If $(\dimv Y)_\infty=0$, then $\mu(Y)=0>\mu(X)$. Thus, $X$ is semistable if and only if it is stable if and only if there is no non-zero subrepresentation $Y$ of $X$ such that $(\dimv Y)_\infty=0$. But giving such a subrepresentation $Y$ is equivalent to giving a subrepresentation $U$ of $M$ such that $U$ is contained in $\Ker f$. We arrive at Definition \ref{ds} of stability for points of $R_{d,n}(Q)$. \hb

\begin{corollary}\label{openpfb} The set $R_{d,n}^s$ is an open subset of $R_{d,n}$. The quotient map $R_{d,n}^s\rightarrow{\cal M}_{d,n}$ is a principal $G_d$-bundle. In particular, it is a smooth morphism.
\end{corollary}

\proof By \cite[Lemma 6.5]{Re}, the action of $G_d$ on $R_{d,n}^s$ is free in the sense of \cite[0.8.iv]{M}, using the above realization of ${\cal M}_{d,n}$. By \cite[Proposition 0.9]{M}, this implies the claimed property.\hb

The technique of extending the quiver used in this section was called ``deframing" by W.~Crawley-Boevey in \cite{CBMW}.

\subsection{Representation-theoretic interpretation}\label{real2}

First we recall the Grassmannians of subrepresentations introduced in \cite{Sc}. Let $X$ be a representation of a quiver $Q$, and let $e\leq d=\dimv X$ be a dimension type. Inside the product of Grassmannians $\prod_{i\in I}{\rm Gr}_{e_i}(X_i)$, we can consider the closed subset ${\rm Gr}_e(X)$ of tuples of subspaces $(U_i)_{i\in I})$ defining a subrepresentation of $X$, that is, fulfilling
$X_\alpha(U_i)\subset U_j$ for all $\alpha:i\rightarrow j$. The subset ${\rm Gr}_e(X)$ is called the Grassmannian of $e$-dimensional subrepresentations of $M$.\\[1ex]
To a point $U=(U_i)_{i\in I}\in{\rm Gr}_e(X)$ we can thus associate a representation of $Q$ of dimension type $e$, again denoted by $U$, in a natural way. Later on, we need to know that this association can also be done on the level of the corresponding varieties:\\[1ex]
The product of Grassmannians $\prod_{i\in I}{\rm Gr}_{e_i}(X_i)$ naturally arises as the quotient of ${\rm IHom}_k(M,X):=\prod_{i\in I}{\rm IHom}_k(M_i,X_i)$ by the group $\prod_{i\in I}{\rm GL}(M_i)$, where ${\rm IHom}_k(M_i,X_i)$ denotes the variety of injective $k$-linear maps from $M_i$ to $X_i$. Denote by ${\rm IHom}_e(X)$ the inverse image of ${\rm Gr}_e(X)$ under this quotient map $\pi$.

\begin{lemma}\label{sigma} There exists a morphism $\sigma:{{\rm IHom}_e(X)}\rightarrow R_e$ such that, for any point $f\in{{\rm IHom}_e(X)}$, the quiver representation $\sigma(f)$ is isomorphic to the quiver representation associated to $\pi(f)$ as above.
\end{lemma}

\proof By definition, a point $f=(f_i)_{i\in I}$ of ${\rm IHom}_k(M,X)$ belongs to ${{\rm IHom}_e(X)}$ if and only if there exist maps $(M_\alpha)_\alpha$ such that $X_\alpha f_i=f_jM_\alpha$ for all arrows $\alpha:i\rightarrow j$ in $Q$. Viewing the $f_i$ as $d_i\times e_i$-matrices, we can cover ${{\rm IHom}_e(X)}$ by open subsets where various tuples of $d_i\times d_i$-minors of these matrices do not vanish. On each of these open subsets, the maps $M_\alpha$, whose existence and uniqueness are guaranteed, can be recovered algebraically via multiplication by inverses of matrices.\hb

For each vertex $i\in I$, denote by $E_i$ the simple representation at $i$, defined by $\dimv E_i=i$. Denote by $P_i$ (resp.~$I_i$) the projective cover (resp.~the injective hull) of $E_i$. These representations can be described explicitly as follows:\\[1ex]
By definition, a path in $Q$ from a vertex $i$ to a vertex $j$ is a sequence of arrows $p=(\alpha_1,\ldots,\alpha_k)$ in $Q$ such that $$i=i_0\stackrel{\alpha_1}{\rightarrow}i_1\stackrel{\alpha_2}{\rightarrow}\ldots\stackrel{\alpha_k}{\rightarrow}i_k=j;$$
we also have one path of length $k=0$ at any vertex $i\in I$. Then the space $(P_i)_j$ (resp.~$(I_i)_j$) has a basis consisting of the paths $p$ from $i$ to $j$ (resp.~from $j$ to $i$). Furthermore, for any representation $M$ of $Q$, there are functorial isomorphisms
$${\rm Hom}_Q(P_i,M)\simeq M_i\mbox{ and }{\rm Hom}_Q(M,I_i)\simeq M_i^*,$$
where $M_i^*$ denotes the $k$-linear dual of the vector space $M_i$.\\[2ex]
Denoting by $I\otimes V$ the injective representation $I\otimes V=\bigoplus_{i\in I}I_i\otimes_k V_i$, we have
$$(I\otimes V)_i\simeq\bigoplus_{j\in I}(I_j)_i\otimes_k V_j\simeq\bigoplus_{j\in I}\bigoplus_{p:i\leadsto j}V_j\simeq \bigoplus_{i\leadsto j}V_j,$$
where $i\leadsto j$ indicates a path from $i$ to $j$.

\begin{definition}\label{defphi}
Given a point $(M,f)\in R_{d,n}$, define a map $\Phi_{(M,f)}=\varphi:M\rightarrow I\otimes V$ recursively as follows:\\[1ex]
If $i$ is a sink in $Q$, then $\varphi_i=f_i:M_i\rightarrow V_i$. For an arbitrary vertex $i\in I$ define
$$\varphi_i=f_i\oplus\bigoplus_{\alpha:i\rightarrow k}\varphi_kM_\alpha:M_i\rightarrow V_i\oplus\bigoplus_{i\rightarrow k}\bigoplus_{k\leadsto j}V_j\simeq\bigoplus_{i\leadsto j}V_j.$$
\end{definition}

\remark We see that we can recover $f$ from $\Phi_{(M,f)}$ by defining $f_i$ as $\varphi_i:M_i\rightarrow\bigoplus_{i\leadsto j}V_j$, followed by the projection to $V_i$.

\begin{lemma}\label{lemma35} The subspace $\Ker\Phi_{(M,f)}$ is the maximal subrepresentation of $M$ contained in $\Ker f$.
\end{lemma}

\proof Suppose $\varphi_im=0$ for some $m\in M_i$. Then $\varphi_jM_\alpha m=0$ for all $\alpha:i\rightarrow j$, thus $M_\alpha m\in\Ker \varphi_j$. This proves that $M_\alpha(\Ker\varphi_i)\subset\Ker\varphi_j$ for all $\alpha:i\rightarrow j$, so that $\Ker\varphi$ is a subrepresentation of $M$.\\[1ex]
Now suppose $U\subset M$ is a subrepresentation contained in the kernel of $f$. We prove by ``downward induction over $Q$" that $U_i\subset\Ker \varphi_i$ for all $i\in I$. If $i\in I$ is a sink in $Q$, then $\varphi_i=f_i$ by definition, thus $U\subset\Ker f_i=\Ker \varphi_i$. For arbitrary $i\in I$, we have $f_i(U)=0$ and $M_\alpha(U_i)\subset U_j\subset\Ker\varphi_j$ by induction, thus $(\varphi_jM_\alpha)(U_i)=0$. By definition of $\varphi_i$, this means $\varphi_i(U_i)=0$. \hb

\begin{corollary}\label{stableinjective} The map $\Phi_{(M,f)}:M\rightarrow I\otimes V$ is injective if and only if the pair $(M,f)$ is stable.
\end{corollary}

\proof The pair $(M,f)$ is stable if and only if there is no non-zero subrepresentation $U$ of $M$ such that $U\subset \Ker f$. But $\Ker\Phi_{(M,f)}$ is the maximal such subrepresentation. \hb

\remark Under the chain of isomorphisms
$${\rm Hom}_Q(M,I\otimes V)\simeq\bigoplus_{i\in I}{\rm Hom}_Q(M,I_i)\otimes_k V_i
\simeq\bigoplus M_i^*\otimes_k V_i\simeq\bigoplus_{i\in I}{\rm Hom}_k(M_i,V_i),$$
the map $f=(f_i)_{i\in I}\in\bigoplus_{i\in I}{\rm Hom}_k(M_i,V_i)$ corresponds to the map $\Phi_{(M,f)}\in{\rm Hom}_Q(M,I\otimes V)$.

\begin{proposition}\label{koko} The framed moduli space ${\cal M}_{d,n}$ is isomorphic to the Grassmannian of subrepresentations ${\rm Gr}_d(I\otimes V)$.
\end{proposition}

\proof Consider the map $\Phi:R_{d,n}\rightarrow{\rm Hom}_k(M,I\otimes V)$ given by $\Phi(M,F)=\Phi_{(M,f)}$. By Corollary \ref{stableinjective}, it maps $R_{d,n}^s$ to the open subset ${\rm IHom}_k(M,I\otimes V)$ of injective linear maps. The map $\Phi$ is also easily seen to be $G_d$-equivariant, for the $G_d$-action on $R_{d,n}$ defined in section \ref{framedmoduli}, and the action on ${\rm Hom}_k(M,I\otimes V)$ given by base change in $M$. By Lemma \ref{lemma35}, the image belongs to ${{\rm IHom}_d(I\otimes V)}$. On the other hand, given a point in ${{\rm IHom}_d(I\otimes V)}$, we can recover both the representation $M$ by Lemma \ref{sigma}, and the map $M\rightarrow V$ by the remark following Definition \ref{defphi}. Thus, we have a $G_d$-equivariant isomorphism $\Phi:R_{d,n}^s\stackrel{\sim}{\rightarrow}{{\rm IHom}_d(I\otimes V)}$. Consequently, the map $\Phi$ descends to an isomorphism of the geometric quotients
$$\overline{\Phi}:{\cal M}_{d,n}\simeq R_{d,n}^s/G_d\rightarrow{{\rm IHom}_d(I\otimes V)}/G_d\simeq{\rm Gr}_d(I\otimes V).$$\hb

\section{Geometry of framed moduli}\label{geometry}

\subsection{Nonemptyness of ${\cal M}_{d,n}$}\label{nonempty}

One of the problems posed by H.~Nakajima in \cite{Na} on framed quiver moduli ${\cal M}_{d,n}$ is to give a criterion for their nonemptyness. We will derive such a criterion using the interpretation of section \ref{real2} and properties of injective coresolutions of representations of quivers.

\begin{lemma}\label{condonni} Given a representation $M$, there exists a map $f:M\rightarrow V$ making the pair $(M,f)\in R_{d,n}$ stable if and only if $n_i\geq\dim{\rm Hom}_Q(E_i,M)$ for all $i\in I$.
\end{lemma}

\proof As in the preceding section, a stable pair $(M,f)$ induces an injection of representations $0\rightarrow M\rightarrow I\otimes V$. Application of the functor ${\rm Hom}_Q(E_i,\_)$ yields an embedding $0\rightarrow{\rm Hom}_Q(E_i,M)\rightarrow{\rm Hom}_Q(E_i,I\otimes V)\simeq V_i$, thus the claimed inequality is neccessary for stability of $(M,f)$.\\[1ex]
We thus have to prove sufficiency. The simple representation $E_i$ admits a projective resolution $\ses{\bigoplus_{i\rightarrow j}P_j}{P_i}{E_i}$. Application of the functor ${\rm Hom}_Q(\_,M)$ gives an exact sequence
$$0\rightarrow{\rm Hom}_Q(E_i,M)\rightarrow\underbrace{{\rm Hom}_Q(P_i,M)}_{\simeq M_i}\rightarrow{\bigoplus_{i\rightarrow j}{\rm Hom}_Q(P_j,M)}{\simeq \bigoplus_{i\rightarrow j}M_j},$$
so that ${\rm Hom}_Q(E_i,M)$ can be identified with $\Ker(\bigoplus_{\alpha:i\rightarrow j}M_\alpha)$. For any $i\in I$, we can choose direct sum decompositions
$$M_i=\Ker(\bigoplus_\alpha M_\alpha)\oplus W_i\mbox{ and }V_i\simeq\Ker(\bigoplus_\alpha M_\alpha)\oplus V_i',$$
since $\dim V_i=n_i\geq\dim{\rm Hom}_Q(E_i,M)=\dim\Ker(\bigoplus_\alpha M_\alpha)$. Define $f_i$ as the projection along $W_i$, followed by the inclusion into $V_i$. Now suppose that $U=(U_i)_i$ is a subrepresentation of $M$ contained in $\Ker(f)=W_i$. We prove $U_i=0$ for all $i\in I$ by descending induction on the partially ordered set $I$. If $i$ is a sink in $Q$, then $W_i=0$, thus $U_i$=0. For arbitrary $i\in I$, we have $M_\alpha(U_i)\subset U_j$ for all $\alpha:i\rightarrow j$. But $U_j=0$ by induction, thus $U_i\subset\Ker\bigoplus_{\alpha:i\rightarrow j}M_\alpha$, which intersects $W_i$ trivially by definition. Thus $U_i=0$. \hb

\begin{corollary}\label{mininj} Any representation $M$ admits an injective coresolution
$$\ses{M}{\bigoplus_{i\in I}I_i\otimes{\rm Hom}_Q(E_i,M)}{\bigoplus_{i\in I}I_i\otimes{\rm Ext}^1_Q(E_i,M)}.$$
\end{corollary}

\proof The above lemma shows the existence of a stable pair $$(M,f:M\rightarrow \bigoplus_{i\in I}{\rm Hom}_Q(E_i,M)),$$ which corresponds to an embedding $0\rightarrow M\rightarrow\bigoplus_{i\in I}I_i\otimes{\rm Hom}_Q(E_i,M)$. Since the path algebra $kQ$ is hereditary, the cokernel is again injective, thus of the form $I\otimes X$ for a graded space $X=(X_i)_i$. Application of ${\rm Hom}_Q(E_i,\_)$ to the resulting exact sequence gives
$$0\rightarrow{\rm Hom}_Q(E_i,M)\rightarrow
\underbrace{{\rm Hom}_Q(E_i,\bigoplus_{j\in I}I_j\otimes{\rm Hom}_Q(E_j,M))}_
{\simeq {\rm Hom}_Q(E_i,M)}\rightarrow$$
$$\rightarrow\underbrace{{\rm Hom}_Q(E_i,\bigoplus_{j\in I}I_j\otimes X_j)}_
{\simeq X_i}\rightarrow{\rm Ext}^1_Q(E_i,M)\rightarrow 0,$$
and thus $X_i\simeq{\rm Ext}^1_Q(E_i,M)$. \hb

\begin{proposition}\label{critnon0} We have $\mathcal{M}_{d,n}\not=\emptyset$ if and only if $n_i\geq \langle i,d\rangle$ for all $i\in I$.
\end{proposition}

\proof That $n_i\geq\langle i,d\rangle$ is a neccessary condition for non-emptyness of ${\cal M}_{d,n}$ follows from Lemma \ref{condonni} since, for any representation $M$ of dimension type $d$, we have $$\langle i,d\rangle=\dim{\rm Hom}_Q(E_i,M)-\dim{\rm Ext}^1_Q(E_i,M).$$
Given $d\in{\bf N}I$, one can clearly define a representation $M$ of dimension type $d$ such that, for all $i\in I$, the rank of the map $\bigoplus_{\alpha:i\rightarrow j}M_\alpha$ is maximal. Then $\dim{\rm Hom}_Q(E_i,M)=\dim\Ker\bigoplus_{\alpha:i\rightarrow j}M_\alpha=\max(0,\langle i,d\rangle)$. This representation thus admits a map $f:M\rightarrow V$ yielding a stable pair $(M,f)$, again by Lemma \ref{condonni}. \hb

\subsection{An embedding of ${\cal M}_{d,n}$}\label{emb}

Recall from section \ref{real2} the embedding $\overline{\Phi}:{\cal M}_{d,n}\rightarrow \prod_{i\in I}{\rm Gr}_{d_i}(\bigoplus_{i\leadsto j}V_j)$ given by $\overline{\Phi}(M,f)=\Imm\Phi_{(M,f)}$.
In this section, we give an explicit description of its image, that is, of the inclusion
$${\rm Gr}_d(I\otimes V)\subset\prod_{i\in I}{\rm Gr}_{d_i}(\bigoplus_{i\leadsto j}V_j).$$

\begin{proposition}\label{propemb} A point $(U_i\subset\bigoplus_{i\leadsto j}V_j)_{i\in I}$ belongs to the image of $\overline{\Phi}$ if and only if for all $i\in I$, we have $U_i\subset V_i\oplus\bigoplus_{\alpha:i\rightarrow j}U_j$.
\end{proposition}

\proof Given $(M,f)\in {\cal M}_{d,n}$, let $\Phi_{(M,f)}=\varphi=(\varphi_i)_{i\in I}$ be the tuple of injective maps constructed in Definition \ref{defphi}, and let $U_i$ be the image of $\varphi_i$. Then
$$U_i=\Imm(f_i\oplus\bigoplus_{\alpha:i\rightarrow j}\varphi_jM_\alpha)\subset\Imm f_i\oplus\bigoplus_{\alpha:i\rightarrow j}\underbrace{\Imm\varphi_jM_\alpha}_{\subset\Imm M_\alpha}\subset V_i\oplus\bigoplus_{\alpha:i\rightarrow j}U_j.$$
Conversely, assume that a tuple of subspaces $U_i\subset \bigoplus_{i\leadsto j}V_j$
of dimensions $\dim U_i=d_i$ is given such that $U_i\subset V_i\oplus\bigoplus_{\alpha:i\rightarrow j}U_j$ for all $i\in I$. We construct a stable pair $(M,f)$ such that $U$ equals the image of $\Imm\Phi_{(M,f)}$. For each vertex $i\in I$, choose a linear map $\varphi_i:M_i\rightarrow\bigoplus_{i\leadsto j}V_j$
with image $U_i$. Define the linear map $f_i$ as the composition $$f_i:M_i\stackrel{\varphi_i}{\rightarrow}\bigoplus_{i\leadsto j}V_j\stackrel{{\rm pr}_{V_i}}{\rightarrow} V_i.$$
For each arrow $\alpha:i\rightarrow j$ in $Q$, consider the composition
$$h_\alpha:M_i\stackrel{\varphi_i}{\rightarrow}\bigoplus_{i\leadsto k}V_k\simeq V_i\oplus\bigoplus_{i\rightarrow l}\bigoplus_{l\leadsto k}V_k\stackrel{{\rm pr}_\alpha}{\rightarrow}\bigoplus_{j\leadsto k}V_k,$$
where ${\rm pr}_\alpha$ denotes the projection map.
By assumption, we have
$$\Imm h_\alpha=\Imm{\rm pr}_\alpha\varphi_i={\rm pr}_\alpha U_i\subset{\rm pr}_\alpha(V_i\oplus\bigoplus_{i\rightarrow k}U_k)\subset U_j=\Imm\varphi_j,$$
thus $h_\alpha$ factors through $\varphi_j$, providing a map $M_\alpha:M_i\rightarrow M_j$ such that  $h_\alpha=M_\alpha\varphi_j$.
From the construction, it is now clear that
$$\varphi_i=f_i\oplus\bigoplus_{\alpha:i\rightarrow j}\varphi_jM_\alpha.$$
Moreover, the maps $\varphi_i$ are injective by definition, thus Lemma \ref{stableinjective} implies that the pair $(M=(M_\alpha)_\alpha,f=(f_i)_{i\in I})$ is stable. \hb

\examples
\begin{enumerate}
\item Let $Q$ be the quiver $i\rightarrow j\rightarrow k$. Then ${\cal M}_{d,n}$ is embedded into
$${\rm Gr}_{d_1}(V_1\oplus V_2\oplus V_3)\times{\rm Gr}_{d_2}(V_2\oplus V_3)\times{\rm Gr}_{d_3}(V_3),$$
and the image consists of those points $(U_1,U_2,U_3)$ of this product of Grassmannians such that
$U_1\subset V_1\oplus U_2$ and $U_2\subset V_2\oplus U_3$.
\item We recall the main example given in \cite{Na}. Let $Q$ be the quiver $1\rightarrow\ldots\rightarrow r$ of type $A_r$. Assume $n_i=0$ for $i\not=r$, and $d_1< d_2<\ldots <d_r<n_r$.
Then ${\cal M}_{d,n}$ is isomorphic to the partial flag variety of flags in $V_r$ with successive steps of dimensions $d_1,d_2,\ldots,d_r,n_r$.
\item More generally, let ${\cal P}$ be a finite poset with a unique maximal element $\infty$. Define a quiver $Q_{\cal P}$ with set of vertices ${\cal P}$, and an arrow $p\rightarrow q$ if and only if $p<q$ in ${\cal P}$ minimally. Define $n_p$ as $0$ if $p\not=\infty$. Then, for arbitrary $d\in{\bf N}{\cal P}$, the variety ${\cal M}_{d,n}(Q_{\cal P})$ parametrizes tuples $(U_p)_{p\in{\cal P}}$ of subspaces of $V_\infty$ such that $U_p\subset U_q$ whenever $p<q$ in ${\cal P}$.
\item Let $Q$ be the quiver $i{\rightarrow\atop\rightarrow}j$. Then ${\cal M}_{d,n}$ parametrises pairs of subspaces $(U_1\subset V_1\oplus V_2^2,U_2\subset V_2)$ such that $U_1\subset V_1\oplus U_2^2$. 
\end{enumerate}

\subsection{${\cal M}_{d,n}$ as an iterated Grassmann bundle}\label{gb}

In this section we describe ${\cal M}_{d,n}$ an an iterated Grassmann bundle. We first recall their definition (see \cite{Fu}). Given a vector bundle $p:E\rightarrow X$ on a variety $X$, we denote by $E_x$ its fibre over a point $x\in X$.

\begin{proposition} Let $X$ be a variety, let $p:E\rightarrow X$ be a vector bundle, and let $d$ be an integer
$0<d<{\rm rk} E$. Then there exists a variety (the Grassmann bundle) ${\rm Gr}_d(E)$, together with a projection $\pi:{\rm Gr}_d(E)\rightarrow X$ and a subbundle $S\subset \pi^*E$ (the tautological bundle), defined by the following universal property:\\[1ex]
For any morphism $f:T\rightarrow X$, the map $g\mapsto g^*S$ induces a functorial bijection between factorizations $g:T\rightarrow {\rm Gr}_d(E)$ of $f$ through $\pi$ and rank $d$-subbundles of $f^*E$. In particular, we have $\pi^{-1}(x)\simeq{\rm Gr}_d(E_x)$
for each $x\in X$.
\end{proposition}

\remark The morphism ${\rm Gr}_d(E)\rightarrow X$ is thus defined as an object in the category of $X$-varieties representing the functor mapping an object $f:T\rightarrow X$ to the set of rank $d$-subbundles of $f^*E$.\\[1ex]
We will use the construction of Grassmann bundles only in the following special situation:

\begin{lemma}\label{cgb} Let $X$ be a variety, and let $E$ be a subbundle of a trivial vector bundle $V\times X$ on $X$. Then the closed subvariety $Y$ of ${\rm Gr}_d(V)\times X$ consisting of pairs $(U,x)$ such that $U\subset E_x$ is isomorphic to ${\rm Gr}_d(E)$, with tautological bundle $S=\{(v,(U,x))\in V\times Y\, :\, v\in U\subset E_x\}$.
\end{lemma}

\proof The universal property defining ${\rm Gr}_d(E)$ can be verified immediately. \hb

To prepare the application of Grassmann bundles to the present setting, we first use the embedding of ${\cal M}_{d,n}$ in $\prod_{i\in I}{\rm Gr}_{d_i}(\bigoplus_{i\leadsto j}V_j)$ to define certain universal bundles ${\cal V}_i$ for $i\in I$ on ${\cal M}_{d,n}$ as follows:

\begin{definition}\label{ub} Denote by ${\cal V}_i$ the subbundle of the trivial bundle $\bigoplus_{i\leadsto j}V_j\times{\cal M}_{d,n}$ on ${\cal M}_{d,n}$ consisting of pairs $(v,(U_i)_i)$ such that $v\in U_i$. For each arrow $\alpha:i\rightarrow j$ in $Q$, define a map ${\cal V}_\alpha:{\cal V}_i\rightarrow{\cal V}_j$ by ${\cal V}_\alpha(v,(U_i)_i)=({\rm pr}_\alpha v,(U_i)_i)$, where ${\rm pr}_\alpha:\bigoplus_{i\leadsto k}V_k\rightarrow\bigoplus_{j\leadsto k}V_k$ denotes the natural projection.
\end{definition}

\remark The discussion of the maps $h_\alpha$ in the proof of Proposition \ref{propemb} shows that the maps ${\cal V}_\alpha$ are indeed well-defined.

\begin{lemma} The tuple $(({\cal V}_i)_{i\in I},({\cal V}_\alpha)_{\alpha})$ is the universal quiver representation on ${\cal M}_{d,n}$, that is, for each point in $(M,f)\in{\cal M}_{d,n}$, the quiver representation induced by the universal one on the fibre over $(M,f)$ is isomorphic to $M$.
\end{lemma}

\proof This follows analogously to the reconstruction of the maps $M_\alpha$ from the maps $h_\alpha$ in the proof of Proposition \ref{propemb}. \hb

Now we come to the central construction in this section:\\[1ex]
Let $i_0$ be a source in the quiver $Q$, and let $\overline{Q}$ be the full subquiver of $Q$ supported on $\overline{I}:=I\setminus\{i_0\}$. Let $\overline{d}$, $\overline{n}$ be the restrictions to $\overline{I}$ of $d$ and $n$, respectively.\\[1ex]
Using Proposition \ref{propemb}, we can realize ${\cal M}_{d,n}(Q)$ as the closed subvariety of ${\rm Gr}_{d_{i_0}}(\bigoplus_{i_0\leadsto j}V_j)\times{\cal M}_{\od,\on}(\oq)$ consisting of pairs $(U,(U_i)_{i\in\oi})$ such that $U\subset V_{i_0}\oplus\bigoplus_{i_0\rightarrow j}U_j$. This induces a projection ${\cal M}_{d,n}(Q)\rightarrow{\cal M}_{\od,\on}(\oq)$. It is also obvious from Definition \ref{ub} that, for all $i\in\oi$, the universal bundles ${\cal V}_{i}$ for $i\in\oi$ of ${\cal  M}_{d,n}(Q)$ and of ${\cal M}_{\od,\on}(\oq)$ are compatible with this realization. Applying Lemma \ref{cgb}, we get immediately:

\begin{lemma}\label{igb} The projection ${\cal M}_{d,n}(Q)\rightarrow{\cal M}_{\od,\on}(\oq)$ is isomorphic to the Grassmann bundle ${\rm Gr}_{d_{i_0}}(V_{i_0}\oplus\bigoplus_{i_0\rightarrow j}{\cal V}_j)$, with tautological bundle ${\cal V}_{i_0}$.
\end{lemma}

Choose an enumeration of $I$ as $\{1,\ldots,s\}$ such that $i<j$ whenever $i\rightarrow j$. Repeated application of the previous lemma yields the main result of this section.

\begin{theorem}\label{chain}
Define a chain of iterated Grassmann bundles
$$M_1\stackrel{p_1}{\rightarrow}M_2\stackrel{p_2}{\rightarrow}\ldots\stackrel{p_{n-1}}{\rightarrow}M_n\stackrel{p_s}{\rightarrow}{\rm pt}$$
by
$$M_i={\rm Gr}_{d_i}(V_i\oplus\bigoplus_{i\rightarrow j}p_{i+1}^*\ldots p_{j-1}^*(S_j)),$$
where $S_i$ denotes the tautological bundle on $M_i$.\\[1ex]
Then ${\cal M}_{d,n}\simeq M_1$, with universal bundles ${\cal V}_i\simeq p_1^*\ldots p_{i-1}^*S_i$ for all $i\in I$.
\end{theorem}

\section{Cohomology}\label{cohomology}

In this section, we will use the realization of Lemma \ref{igb} of the framed moduli ${\cal M}_{d,n}$ to determine its Chow ring. We collect some results on Chow rings of Grassmann bundles from \cite{Fu}.\\[1ex]
Slightly different notation from the one in \cite{Fu} will be convenient, which we introduce now. Given a smooth irreducible variety $X$, let $A^*(X)=\oplus_iA^i(X)$ be its Chow ring, where $A^i(X)$ denotes the group of codimension $i$ cycles modulo rational equivalence, and the ring structure is given by the intersection product. Given a vector bundle $E$ on $X$, its Chern classes $c_i(E)$ for $i\geq 0$ will be considered as elements of $A^i(X)$ (they are considered as degree $i$ endomorphisms of $A^*(X)$ in \cite{Fu}; applying these to the fundamental class $[X]\in A^0(X)$ yields the elements considered here). We also consider the formal Chern polynomial $c_t(E):=\sum_{i=0}^\infty c_i(X)t^i\in A^*(X)[[t]]$.\\[1ex]
Given any formal power series $P=\sum_{i=0}^\infty c_it^i\in A[[t]]$ with coefficients in a commutative ring $A$ and a partition $\lambda=(\lambda_1\geq\ldots\geq\lambda_d\geq 0)$, we define the formal Schur polynomial $\Delta_\lambda(P)\in A$ as the determinant of the $d\times d$-matrix $(c_{\lambda_i+j-i})_{i,j}$.

\begin{theorem}[\cite{Fu}] Let $X$ be a smooth irreducible variety, let $E$ be a vector bundle on $X$, and let $0<d<{\rm rk} E$ be an integer. We consider the corresponding Grassmann bundle $\pi:{\rm Gr}_d(E)\rightarrow X$ with tautological bundle $S$.
\begin{enumerate}
\item For all $i\in{\bf N}$, there is an isomorphism $A^i({\rm Gr}_d(E))\simeq\oplus_\lambda A^{i-|\lambda|}(X)$, where the direct sum runs over partitions $\lambda=({\rm rk} E-d\geq\lambda_1\geq\ldots\geq\lambda_d\geq 0)$, and $|\lambda|=\sum_k\lambda_k$. Each element of $A^i({\rm Gr}_d(E))$ has a unique expression in the form $\sum_\lambda\Delta_\lambda(c_t(S)^{-1})\cdot \pi^*\alpha_\lambda$ for elements $\alpha_\lambda\in A^{i-|\lambda|}(X)$, the sum being as above.
\item The ring $A^*({\rm Gr}_d(E))$ is generated by $A^*(X)$, together with the Chern classes $c_1(S),\ldots,c_d(S)$, with the following defining relation: the formal power series $c_t(S)^{-1}$ is a polynomial of degree at most ${\rm rk} E-d$.
\item For the base field $k={\bf C}$, the cycle map $A^*({\rm Gr}_d(E))\rightarrow H^*({\rm Gr}_d(E))$ to singular cohomology is an isomorphism (in particular, there is no odd cohomology) if and only if the cycle map $A^*(X)\rightarrow H^*(X)$ is so.
\item The Poincar\'e polynomial $P_{{\rm Gr}_d(E)}(q):=\sum_{i=0}^\infty\dim A^i(X)q^i$ is given by $P_X(q)\cdot\left[{{\rm rk E}\atop d}\right]$, where $\left[{n\atop d}\right]=\prod_{k=1}^d\frac{q^{n-d+k}-1}{q^k-1}$ denotes the $q$-binomial coefficient.
\end{enumerate}
\end{theorem}

\proof The first part is precisely \cite[Proposition 14.6.5]{Fu} translated to our conventions. Note that the Schur polynomials there are taken in a certain formal power series $c_t(Q-\pi^*E)$ in terms of the sequence of vector bundles $0\rightarrow S\rightarrow \pi^*E\rightarrow Q\rightarrow 0$, which identifies with $c_t(S)^{-1}$ by the Whitney sum formula \cite[Theorem 3.2]{Fu}. The second part is a reformulation of \cite[Example 14.6.6]{Fu}, where the ring structure is described by generators $c_1(S),\ldots,c_d(S)$ and $c_1(Q),\ldots,c_{{\rm rk E}-d}(Q)$ (for the bundle $Q$ defined above) subject to the relation $c_t(S)\cdot c_t(Q)=c_t(\pi^*(E))$. Thus, as a formal power series, $c_t(Q)$ is given by $c_t(\pi^*(E))/c_t(S)$, the only restriction being that this is actually a polynomial of degree $\leq {\rm rk E}-d$. The third part is \cite[Example 19.1.11 d)]{Fu}, together with Poincare duality. The fourth part follows from the first by some standard combinatorics.\hb

Repeated application of this theorem to the chain of Grassmann bundles in Theorem \ref{chain} easily yields the following description of the Chow ring of ${\cal M}_{d,n}$:

\begin{theorem}\label{coho} For all ${\cal M}_{d,n}\not=0$, the following holds:
\begin{enumerate}
\item The Chow ring $A^\bullet({\cal M}_{d,n})$ has a linear basis
$$\{\prod_{i\in I}\Delta_{\lambda^i}(c_t({\cal V}_i)^{-1}),$$
indexed by tuples of partitions $(\lambda^i)_{i\in I}$ such that
$$n_i-\langle i,d\rangle\geq\lambda_1^i\geq\ldots\geq \lambda_{d_i}^i\geq 0.$$
\item As a ring, $A^\bullet({\cal M}_{d,n})$ is generated by the Chern classes of the universal bundles ${\cal V}_i$ for $i\in I$, together with the following defining relations:\\[1ex]
For each $i\in I$, the formal power series
$$\frac{\prod_{i\rightarrow j}c_t({\cal V}_j)}{c_t({\cal V}_i)}$$
is a polynomial of degree at most $\leq n_i-\langle i,d\rangle$.
\item Assume $k={\bf C}$. Then the odd cohomology $H^{2\bullet+1}({\cal M}_{d,n})$ vanishes, and the even cohomology $H^{2\bullet}({\cal M}_{d,n})$ is isomorphic to the Chow ring $A^\bullet({\cal M}_{d,n})$.
\item The Poincar\'e polynomial of ${\cal M}_{d,n}$ is given by
$$\sum_{k=0}^\infty\dim A^{k}({\cal M}_{d,n})q^k=\prod_{i\in I}\left[{{n_i+\sum_{i\rightarrow j}d_j}\atop{d_i}}\right].$$

\end{enumerate}
\end{theorem}

\remark The third part of the theorem is also implied by a theorem of \cite{KW}, where this property is proved in the generality of quiver moduli for which stability and semistability coincide. There it is also proved that, in this case, the Chern classes of the universal bundles generate the Chow ring (see also \cite{ES,Na}).\\[1ex]
Note that the description of a basis of the Chow ring in terms of Schur polynomials allows application of all classical formulas of Schubert calculus, like for example the Pieri rule, Giambelli's formula and the Littlewood-Richardson rule (see \cite{Fu}).

\section{Orbit structure}\label{orb}

In this section, we define a group action on the framed moduli ${\cal M}_{d,n}$ and relate it to the action of the group $G_d$ on the representation variety $R_d$.\\[1ex]
Recall from the proof of Proposition \ref{koko} the $G_d$-equivariant isomorphism
$\Phi: R_{d,n}^s\stackrel{\sim}{\rightarrow}{\rm IHom}_d(I\otimes V)$. On the right hand side, we have a natural action of ${\rm Aut}_Q(I\otimes V)$, which translates into an action on $R_{d,n}^s$. This latter action obviously commutes with the $G_d$-action, thus it induces an action of ${\rm Aut}_Q(I\otimes V)$ on the quotient ${\cal M}_{d,n}$.\\[1ex]
To identify this action in terms of the embedding ${\cal M}_{d,n}\subset\prod_{i\in I}{\rm Gr}_{d_i}(\bigoplus_{i\leadsto j}V_j)$, we consider
$${\rm End}_Q(I\otimes V)\simeq\bigoplus_{i,j}{\rm Hom}_Q(I_i\otimes V_i,I_j\otimes V_j)\simeq\bigoplus_{i\leadsto j}{\rm Hom}_k(V_i,V_j)$$
as a subspace of
$$\bigoplus_{i\in I}{\rm End}_k((I\otimes V)_i)\simeq{\rm End}_k(\bigoplus_{i\leadsto j}V_j).$$
Since $Q$ has no oriented cycles, the group ${\rm Aut}_Q(I\otimes V)$ identifies with the set of tuples $(f_\omega)_\omega\in\bigoplus_{\omega:i\leadsto j}{\rm Hom}_k(V_i,V_j)$ such that all components $f_{\omega:i\leadsto i}$ corresponding to the empty paths are invertible.\\[1ex]
We thus arrive at an embedding of
$$A_n:=\{(f_\omega)_\omega\in\bigoplus_{\omega:i\leadsto j}{\rm Hom}_k(V_i,V_j)\, ;\, \mbox{all }f_{\omega:i\leadsto i}\mbox{ invertible}\}$$
into $\prod_{i\in I}{\rm GL}_k(\bigoplus_{i\leadsto j}V_j)$, which proves:

\begin{lemma} The canonical action of the group $\prod_{i\in I}{\rm GL}_k(\bigoplus_{i\leadsto j}V_j)$ on the variety $\prod_{i\in I}{\rm Gr}_{d_i}(\bigoplus_{i\leadsto j}V_j)$ induces an action of the subgroup $A_n\simeq{\rm Aut}_Q(I\otimes V)$ on ${\cal M}_{d,n}$.
\end{lemma}

\example For the quiver $Q=i\rightarrow j\rightarrow k$, the group $A_n$ is the subgroup of
$${\rm GL}(V_1\oplus V_2\oplus V_3)\times{\rm GL}(V_2\oplus V_3)\times{\rm GL}(V_3)$$
given by block matrices of the form
$$(\left[\begin{array}{ccc}a&&\\ b&d&\\ c&e&f\end{array}\right],\left[\begin{array}{cc}d&\\ e&f\end{array}\right],[f])$$
such that the matrices $a$, $d$ and $f$ are invertible.\\[1ex]
To analyse this group action on ${\cal M}_{d,n}$, we first record two simple representation-theoretic facts.

\begin{lemma} Given a representation $M$ and two embeddings $i,j$ of $M$ into an injective representation $I$, any automorphism $\varphi$ of $M$ extends to an automorphism $\psi$ of $I$, such that $j\varphi=\psi i$.
\end{lemma}

\proof The image of $i$ being contained in an injective hull $I'$ of $M$, we can write $I=I'\oplus I''$ and $i=\left[{i'\atop  0}\right]$. The map $j\varphi$ factors though $i$, yielding maps $\psi':I'\rightarrow I$ and $\psi'':I''\rightarrow I$ such that $\psi'i'=j\varphi$. Since $I'$ is an injective hull of $M$, injectivity of $j\varphi$ implies injectivity of $\psi'$. Without loss of generality we can assume $\psi''$ to be injective, too. The map $\psi=[\psi',\psi'']$ is the claimed extension.\hb

\begin{lemma}\label{extend} Two subrepresentations of an injective representation $I$ are conjugate under the action of ${\rm Aut}_Q(I)$ if and only if they are isomorphic. The stabilizer ${\rm Iso}_{{\rm Aut}_Q(I)}(U)$ is the extension of the automorphism group ${\rm Aut}_Q(U)$ by a unipotent group.
\end{lemma}

\proof If two subrepresentations are conjugate under ${\rm Aut}_Q(I)$, they are clearly isomorphic. The converse follows from the above lemma. To analyse the stabilizer, choose an embedding $i:M\rightarrow I$ and consider the group $G:=\{(\varphi,\psi)\in{\rm Aut}_Q(M)\times{\rm Aut}_Q(I)\, :\, i\varphi=\psi i\}$. Projection to the second component yields an isomorphism between $G$ and ${\rm Stab}_{{\rm Aut}_Q(I)}(U)$, since the automorphism $\varphi$ is uniquely determined by $\psi$, if it exists. Again by the above lemma, projection to the first component yields a surjective map to ${\rm Aut}_Q(M)$, whose kernel consists of the automorphisms $\varphi$ of $M$ such that $\varphi i=i$. This condition is equivalent to $\varphi$ being contained in $1+{\rm Hom}_Q(I/M,I)p$, where $p:I\rightarrow I/M$ denotes the canonical projection. This last group is clearly unipotent.\hb

Denote by $R_{d}^{(n)}$ the image of the natural projection $p:R_{d,n}^s\rightarrow R_d$. By Lemma \ref{condonni}, $R_d^{(n)}$ consists of the representations $M$ such that $\dim{\rm Hom}_Q(E_i,M)\leq n_i$ for all $i\in I$. Since, by Corollary \ref{openpfb}, the map $p$ can be factored into an open embedding followed by a projection along an affine space, it is smooth. Furthermore, the quotient map $\pi:R_{d,n}^s\rightarrow {\cal M}_{d,n}$ is smooth by Corollary \ref{openpfb}.\\[1ex]
Thus, we arrive at a diagram of smooth morphisms
$$R_d^{(n)}\stackrel{p}{\leftarrow}R_{d,n}^s\stackrel{\pi}{\rightarrow}{\cal M}_{d,n}.$$ 
The map $p$ is $G_d$-equivariant and $A_n$-invariant, whereas the map $\pi$ is $A_n$-equivariant and a $G_d$-quotient. We can therefore conclude:

\begin{theorem}\label{bij} The map $X\mapsto \pi p^{-1}X$ gives a bijection between $A_n$-stable subvarieties of ${\cal M}_{d,n}$ and $G_d$-stable subvarieties of $R_d^{(n)}$, such that inclusions, closures, irreducibility and types of singularities are preserved.
\end{theorem}

\proof The map $X\mapsto \pi p^{-1}X$ maps a single $G_d$-orbit in $R_d$ to a single $A_n$-orbit in ${\cal M}_{d,n}$ by Lemma \ref{extend}. By definition of the subset $R_d^{(n)}$, all $A_n$-orbits in ${\cal M}_{d,n}$ arise in this way. This proves bijectivity of the correspondence. All the geometric properties follow from the above mentioned properties of the morphisms $p$ and $\pi$.\hb

\remark This result is also obtained, in the generality of representations of finite dimensional algebras, in \cite{Bo}.\\[1ex]
The above theorem allows us to derive several results on the $A_n$-orbit structure of ${\cal M}_{d,n}$ from corresponding results on the action of $G_d$ on $R_d$. We just give two example applications.\\[1ex]
Applying Gabriel's Theorem, we get immediately:

\begin{corollary} If $Q$ is of Dynkin type, that is, the underlying unoriented graph is a disjoint union of Dynkin diagrams of type $A_l$, $D_l$, $E_6$, $E_7$, $E_8$, then $A_n$ acts on ${\cal M}_{d,n}$ with finitely many orbits.
\end{corollary}

\begin{corollary} If the quiver $Q$ is of type $A$ or $D$, all closures of $A_n$-orbits in ${\cal M}_{d,n}$ are normal and Cohen-Macaulay.
\end{corollary}

\proof This follows from the corresponding properties of closures of $G_d$-orbits in $R_d$, proven in \cite{Zw1}, \cite{Zw2}.\hb

In view of the description of the Chow ring of ${\cal M}_{d,n}$ of Theorem \ref{coho}, one can ask the following question:\\[1ex]
What is the class $[\overline{\cal O}]\in A^*({\cal M}_{d,n})$ of the closure of an $A_n$-orbit in ${\cal M}_{d,n}$, given either in terms of the linear basis of products of Schur polynomials, or in terms of the Chern classes of universal bundles?\\[2ex]
\example If $Q$ is the equioriented quiver $0\rightarrow 1\rightarrow\ldots\rightarrow n$ of type $A_{n+1}$,
we can compute, at least in principle, the classes of orbit closures in the Chow ring, applying the formula for degeneracy loci of \cite{BF}. This can be done as follows:\\[1ex]
Using Theorem \ref{bij} and the description of the closures of $G_d$-orbits in $R_d$ of \cite{AD}, we see that the closures of $A_n$-orbits are precisely of the form
$$\{((M,f)\in{\cal M}_{d,n}\, :\, {\rm rk}(M_{ij})\leq r_{ij}\mbox{ for all }0\leq i<j\leq n\},$$
where $M_{ij}$ denotes the composite map $M_i\rightarrow\ldots\rightarrow M_j$, and the tupel $r=r_{ij}$ is as in \cite[(1.2)]{BF}.
By the description of the universal quiver representation ${\cal V}$ on ${\cal M}_{d,n}$ of Definition \ref{ub}, we see that this orbit closure is thus precisely the degeneracy locus $\Omega_r(({\cal V})$ of \cite[(1.1)]{BF}. Thus, by the main result of \cite{BF}, its class in the Chow ring can be written in the form
$$\sum_\lambda c_\lambda(r)\Delta_{\lambda^1}(c_t({\cal V}_1)/c_t({\cal V}_0))\ldots\Delta_{\lambda^n}(c_t({\cal V}_n)/c_t({\cal V}_{n-1})),$$
the sum running over tuples $\lambda^1,\ldots,\lambda^n$ of partitions, and the $c_\lambda(r)$ being certain integer coefficients. Furthermore, by \cite[Lemma 4]{BF}, we have for any two bundles $E$, $F$ on a variety $X$ the formula
$$\Delta_\lambda(c_t(E)/c_t(F))=\sum_{\mu,\nu}c_{\mu,\nu}^\lambda\Delta_{\mu}(c_t(E)^{-1})\Delta_\nu(c_{-t}(F)^{-1}),$$
the $c_{\mu,\nu}^\lambda$ denoting Littlewood-Richardson coefficients. These two formulas allow us to express the class of the orbit closure in terms of the basis exhibited in Theorem \ref{coho}.

\section{Realizations of quantum groups}\label{quantum}

We first recall the quantum group construction of ${\cal U}_q(\mathfrak{gl}_n)$ of \cite{BLM}. Let $k$ be a finite field, and let $v\in{\bf C}$ be a square root of the cardinality of $k$. Let ${\cal U}_v(\mathfrak{gl}_n)$ be the quantized enveloping algebra of $\mathfrak{gl}_n$, specialized at the quantum parameter $v\in{\bf C}$ (see \cite{BLM}).\\[1ex]
Let $V$ be a $k$-vector space of dimension $N\in{\bf N}$. Let ${\cal F}(V)$ be the variety of $n$-step flags
$$0=V_0\subset V_1\subset\ldots\subset V_d=V$$
in $V$, which is a non-connected projective $k$-variety. We have an obvious action of ${\rm GL}(V)$ on ${\cal F}(V)$, thus a diagonal action on ${\cal F}(V)\times{\cal F}(V)$. Denote by
$$A={\bf C}^{{\rm GL}(V)}({\cal F}(V)\times{\cal F}(V))$$
the space of (arbitrary) ${\rm GL}(V)$-invariant functions from ${\cal F}(V)\times{\cal F}(V)$ to ${\bf C}$. The space $A$ carries an associative multiplication given by convolution:
$$(f\cdot g)(F_1,F_3)=\sum_{F_2\in{\cal F}(V)}f(F_1,F_2)g(F_2,F_3)$$
for functions $f,g\in A$, and for all $F_1,F_3\in{\cal F}(V)$.

\begin{theorem}[\cite{BLM}] There exists a homomorphism of ${\bf C}$-algebras
$${\cal U}_v(\mathfrak{gl}_n)\rightarrow A(V).$$
\end{theorem}

The aim of this section is to develop a partial generalization of this result. More precisely, let $C=(\langle i,j\rangle+\langle j,i\rangle)_{i,j\in I}$ be the matrix representing the symmetrization of the Euler form. This being a symmetric generalized Cartan matrix, we have an associated Kac-Moody Lie algebra $\mathfrak{g}$ with triangular decomposition $\mathfrak{g}=\mathfrak{n}^-\oplus\mathfrak{h}\oplus\mathfrak{n}^+$. We denote by ${\cal U}_v(\mathfrak{n}^+)$ the quantized enveloping algebra of $\mathfrak{n}^+$ with quantum parameter $v$, which is a complex algebra. Our aim is to realize certain quotients of a modified version $\dot{\cal U}_v(\mathfrak{b})$ of ${\cal U}_v(\mathfrak{n}^+)$ (see below for the definition) in terms of a convolution algebra related to the varieties ${\cal M}_{d,n}$.\\[2ex]
We continue to denote by $k$ a finite field with $v^2$ elements. We can easily repeat the definition of $R_d$, ${\cal M}_{d,n}$ etc. over $k$ {\it as sets}. It is then easy to see that all set-theoretic results, like the description of the embedding of ${\cal M}_{d,n}$ into $\prod_{i\in I}{\rm Gr}_{d_i}(\bigoplus_{i\leadsto j}V_j)$, the existence of the group action $A_n:{\cal M}_{d,n}$ etc. continue to hold. So we will use all these notions and results freely in this section.

\subsection{Review of Hall algebras}\label{hall}

We recall the Hall algebra construction of C.~M.~Ringel \cite{Ri}. Let ${\bf C}^{G_d}(R_d)$ be the space of $G_d$-invariant functions from $R_d$ to ${\bf C}$. Define
$$H=H_v(Q)=\bigoplus_{d\in{\bf N}I}{\bf C}^{G_d}(R_d).$$
We define a ${\bf N}I$-graded convolution on $H$ as follows: let $f\in{\bf C}^{G_d}(R_d)$, $g\in{\bf C}^{G_e}(R_e)$ be two functions in $H$. Define their convolution $f\cdot g\in{\bf C}^{G_{d+e}}(R_{d+e})$ by 
$$(f\cdot g)(X)=v^{\langle e,d\rangle}\sum_{U\subset X}f(U)g(X/U),$$
where the sum runs over all subrepresentations $U$ of the representation $X$. Note that this is well-defined, since the sum is obviously finite, and the value of $f$ on $U$ (resp.~of $g$ on $X/U$) is well-defined since it only depends on the isomorphism class of $U$ (resp.~$X/U$).\\[1ex]
Let $E_i$ be the function on $R_i$ sending the unique point to $1$. Define $C=C_v(Q)$ as the ${\bf C}$-subalgebra of $H$ generated by the functions $E_i$ for $i\in I$. We have

\begin{theorem}[\cite{Gr,Ri}] The map $\eta$ sending the Chevalley generator $E_i$ of ${\cal U}_v^+(\mathfrak{g})$ to the function $E_i$ induces an isomorphism
$$\eta:{\cal U}_v^+(\mathfrak{g})\stackrel{\sim}{\rightarrow}C.$$
Moreover, if $Q$ is of Dynkin type, then $C=H$.
\end{theorem}

For the following, we need a modified version of the Hall algebra, reminiscent of the modified quantized enveloping algebra of \cite{Lu}.

\begin{definition} Define a modified Hall algebra $\dot{H}=\dot{H}_v(Q)$ by adjoining to the algebra $H$ pairwise orthogonal idempotents $1_d$ for $d\in{\bf Z}I$ fulfilling the relation $f1_d=1_{d+|f|}$ for all $d\in{\bf Z}I$ and all homogeneous elements $f\in H$. Define $\dot{C}=\dot{C}_v(Q)$ as the subalgebra of $\dot{H}$ generated by the functions $E_i$ and all idempotents $1_d$ for $d\in{\bf Z}I$.
\end{definition}

Note that, by the defining relations of $\dot{H}$, any element in $\dot{H}$ can be written uniquely in the form $\sum_{d\in{\bf Z}I}1_df_d$ for elements $f_d\in H$, almost all of which are zero. The multiplication in $\dot{H}$ is then given by the rule $1_df\cdot 1_eg=\delta_{d,e+|f|}1_d(f\cdot g)$ for homogeneous elements $f,g\in H$. More precisely, the Hall algebra $H$ has a basis consisting of the characteristic functions $\chi_{[M]}$ of the $G_d$-orbits in the spaces $R_d$, which correspond bijectively to the isomorphism classes $[M]$ of representations of $Q$. Thus, we have a corresponding basis $1_d\cdot\chi_{[M]}$ of $\dot{H}$.\\[1ex]
In the same way as for $H$, we define $\dot{\cal U}_v(\mathfrak{b})$ by adjoining orthogonal idempotents to ${\cal U}_v(\mathfrak{n^+})$ (see \cite{Lu}). Then we have:

\begin{corollary}\label{raute} The map $\eta$ above extends to an isomorphism $\eta:\dot{\cal U}_v(\mathfrak{b})\stackrel{\sim}{\rightarrow}\dot{C}$. Moreover, if $Q$ is of Dynkin type, we have $\dot{C}=\dot{H}$.
\end{corollary}

\subsection{Quantum group construction from ${\cal M}_{d,n}$}\label{qgc}

\begin{definition}
\begin{enumerate}
\item Denote by ${\cal M}_{d,n}^{\rm proj}$ the set of all pairs $(M,f)$ in ${\cal M}_{d,n}$, such that $M$ is a projective representation of $Q$.
\item Denote by ${\cal M}_n^{\rm proj}$ the disjoint union of the ${\cal M}_{d,n}^{\rm proj}$, for all $d\in{\bf N}I$.
\item Finally, denote by ${\cal X}_n$ the subvariety of ${\cal M}_n^{\rm proj}\times {\cal M}_n^{\rm proj}$ of pairs of tuples of subspaces $((U_i)_i,(U'_i)_i)$ such that $U_i$ contains $U'_i$ for all $i\in I$.
\end{enumerate}
\end{definition}

\remark Since the projective representations in $R_d$ form a single $G_d$-orbit (if any), the subset ${\cal M}_{d,n}^{\rm proj}$ is either empty or a single $A_n$-orbit by the correspondence of Theorem \ref{bij}. By the criterion \ref{critnon0} for non-emptyness of ${\cal M}_{d,n}$, the set ${\cal M}_n$ is a finite union of the ${\cal M}_{d,n}^{\rm proj}$. Therefore, $A_n$ acts on ${\cal M}_{n}$ with finitely many orbits.\\[1ex]
Representation-theoretically, a point in ${\cal X}_n$ can be viewed as a diagram
$$P'\rightarrow P\rightarrow I\otimes V$$
of embeddings of representations of $Q$, where $P$ and $P'$ are both projective.\\[2ex]
We first define a map $\xi$ from the set ${\cal X}_n/A_n$ of $A_n$-orbits in ${\cal X}_n$ to the disjoint union
$\dot{\bigcup}_{d\in{\bf N}I}R_d/G_d$ of all $G_d$-orbits in all $R_d$ as follows:

\begin{definition} Given a point $(U_*,U'_*)=((U_i)_i,(U_i')_i)$ in ${\cal X}_n$,
let $\xi(U_*,U'_*)=U_*/U'_*$ be the representation of $Q$ defined by the following tuple of $k$-vector spaces
and $k$-linear maps:
$$U_*/U'_*=((U_i/U'_i)_{i\in I},(\overline{{\rm pr}_\alpha}:U_i/U'_i\rightarrow U_j/U'_j)_{\alpha:i\rightarrow j}).$$
\end{definition}

We endow $C:={\bf C}^{A_n}({\cal X}_n)$ with an algebra structure via convolution of functions:
$$(f\cdot g)(U_*,U''_*)=\sum_{U'_*}f(U_*,U'_*)g(U'_*,U''_*).$$
This algebra is ${\bf N}I$-graded by defining the homogeneous component $C_d$ as the space of functions supported on $(\bigcup_{e\in{\bf N}I}{\cal M}_{d+e,n}\times{\cal M}_e)\cap {\cal X}_n$.\\[1ex]
Furthermore, we define the map $\dimv$ on ${\cal M}_{n}$ by $\dimv(M,f)=\dimv M$.

\begin{proposition} The map $\xi$ induces a morphism $\eta:H\rightarrow C$ of ${\bf N}I$-graded ${\bf C}$-algebras
defined by
$$\eta(f)(U_*,U'_*)=v^{\langle \dimv U_*-\dimv U'_*,\dimv U'_*\rangle}f(U_*/U'_*).$$
\end{proposition}

\proof Well-definedness of the map $\eta$ follows from the above definition of the map $\xi$.
The map $\eta$ preserves the ${\bf N}I$-gradings by their definition.
It remains to show compatibility with the respective convolution structures.
So let $f$ and $g$ be elements of $H$ of degree $d$ and $e$, respectively, and let $(U_*,U''_*)$ be a point in ${\cal X}_n$.
For the calculation of $\eta(f)(U_*,U''_*)$, we can assume without loss of generality that $\dimv U_*=\dimv U''_*+d+e$; denote $c=\dimv U''_*$. By the definitions, we have
$$\eta(f\cdot g)(U_*,U''_*)=v^{\langle d+e,c\rangle}(f\cdot g)(M)=
v^{\langle d+e,c\rangle+\langle e,d\rangle}\sum_{Z}f(Z)g((U_*/U''_*)/Z),$$
where the sum runs over all subrepresentations $Z$ of $U_*/U''_*$ of dimension type $\dimv M=d$. By definition of the representation $U_*/U'_*$, such subrepresentations are in bijection with points $U'_*$ in ${\cal M}^{\rm proj}_{c+d}$ such that $U_*\supset U'_*\supset U''_*$ via the identification $Z=U'_*/U''_*$. Thus, we can rewrite this last sum as
$$v^{\langle e,c+d\rangle+\langle d,c\rangle}\sum_{U_*\supset U'_*\supset U''_*}f(U_*/U'_*)g(U'_*/U''_*)=\sum_{U'_*}\eta(f)(U_*,U'_*)\eta(g)(U'_*,U''_*),$$
which equals $(\eta(f)\cdot \eta(g))(U_*,U''_*)$, as stated.\hb

For each $d\in{\bf Z}I$, define a function $\Delta_d$ on ${\cal X}_n$ by $\Delta_d(U_*,U'_*)=0$ if $U_*\not=U'_*$ or $\dimv U_*\not=d$, and $\Delta_d(U_*,U_*)=1$ if $\dimv U_*=d$.

\begin{lemma} The map $\eta$ extends to an algebra morphism, again denoted by $\eta$, from $\dot{H}$ to ${\bf C}^{A_n}({\cal X}_n)$, mapping $1_d$ to the function $\Delta_d$.
\end{lemma}

\proof Obviously, the functions $\Delta_d$ are pairwise orthogonal idempotents in ${\bf C}^{A_n}({\cal X}_n)$. The relation $f1_d=1_{d+e}f$ is immediately verified.\hb

To understand the image of the map $\eta$, we have to analyse the set of $A_n$-orbits in ${\cal X}_n$.

\begin{lemma}\label{charo} The $A_n$-orbit of a point $(U_*,U'_*)$ in ${\cal X}_n$ is uniquely determined by the isomorphism class of the representation $U_*/U'_*$, together with the dimension type $\dimv U_*\in{\bf N}I$.
\end{lemma}

\proof Clearly, both the isomorphism class $U_*/U'_*$ and the dimension type $\dimv U_*$ are constant along the $A_n$-orbit of $(U_*,U'_*)$. Conversely, assume that for points $(^1U_*,^1U'_*)$ and $(^2U_*,^2U'_*)$ in ${\cal X}_n$, the representations $
~^1U_*/^1U'_*$ and $~^2U_*/^2U'_*$ are isomorphic, and $\dimv^1U_*=\dimv^2U_*$. We have to prove that these two points are conjugate under $A_n$. We choose diagrams
$$P_1\stackrel{\alpha_1}{\rightarrow}P'_1\stackrel{i_1}{\rightarrow}I\otimes V\mbox{ and }P_2\stackrel{\alpha_2}{\rightarrow}P'_2\stackrel{i_2}{\rightarrow}I\otimes V,$$
respectively, representing these two points. By assumption, we have $\dimv P_1=\dimv P_2$ and $P_1/P'_1\simeq P_2/P'_2=:M$. We thus have exact sequences $$\ses{P'_1}{P_1}{M}\mbox{ and }\ses{P'_2}{P_2}{M}.$$
Using the obvious dual version of Lemma \ref{extend}, we can extend the identity map on $M$ to an isomorphism $g:P_1\rightarrow P_2$; this in turn yields an isomorphism $f:P'_1\rightarrow P'_2$ such that $g\alpha_1=\alpha_2f$. Applying Lemma \ref{extend} again, we get an automorphism $\psi$ of $I\otimes V$ extending $g$, that is, $\psi i_1=i_2g$. With the calculation $$\psi i_1\alpha_1=i_2g\alpha_1=i_2\alpha_2 f$$ we see that the map $\psi$ also extends the map $f$ along the embeddings $i_1\alpha_1$ and $i_2\alpha_2$, respectively. But this means precisely that $\psi\in{\rm Aut}_Q(I\otimes V)=A_n$ conjugates the two given points in ${\cal X}_n$. The following commutative diagrams, the left one having exact rows, illustrate this construction:
$$\begin{array}{ccccccccccccccc}
0&\rightarrow&P'_1&\stackrel{\alpha_1}{\rightarrow}&P_1&\rightarrow&M&\rightarrow&0&&
P'_1&\stackrel{\alpha_1}{\rightarrow}&P_1&\stackrel{i_1}{\rightarrow}&I\otimes V\\
&&f\downarrow&&g\downarrow&&||&&&&
f\downarrow&&g\downarrow&&\psi\downarrow\\
0&\rightarrow&P'_2&\stackrel{\alpha_2}{\rightarrow}&P_2&\rightarrow&M&\rightarrow&0&&
P'_2&\stackrel{\alpha_2}{\rightarrow}&P_2&\stackrel{i_2}{\rightarrow}&I\otimes V
\end{array}$$\hb

This characterization of the $A_n$-orbits in ${\cal X}_N$ leads us to define, for all $d\in{\bf Z}I$ and all isomorphism classes $[M]$ of representations of $Q$, the set
$${\cal O}_{d,[M]}:=\{(U_*,U'_*)\in{\cal X}_n\, :\, \dimv U_*=d,\;\; U_*/U'_*\simeq M\}.$$

\begin{proposition} The set ${\cal O}_{d,[M]}$ is non-empty if and only if $\langle i,d\rangle \leq n_i$ and $\dim{\rm Hom}_Q(M,E_i)\leq\langle d,i\rangle$ for all $i\in I$. In this case, it is a single $A_n$-orbit in ${\cal X}_n$, and all such orbits arise in this way.
\end{proposition}

\proof By Lemma \ref{charo}, we only have to prove the criterion for non-emptyness of ${\cal O}_{d,[M]}$. So suppose that ${\cal O}_{d,[M]}\not=\emptyset$. Any point $(U_*,U'_*)$ in this set then provides projective representations $P$ and $P'$, such that there exists an exact sequence $\ses{P'}{P}{M}$, and such that $\dimv P=d$. By Proposition \ref{critnon0}, we neccessarily have $\langle i,d\rangle\leq n_i$ for all $i\in I$. Furthermore, we have $\dim{\rm Hom}_Q(M,E_i)\leq\dim{\rm Hom}_Q(P,E_i)=\langle d,i\rangle$. This proves neccessity of the numerical conditions.\\[1ex]
To prove the converse, assume that the above two sets of conditions on $d$ and $M$ are satisfied. By the obvious dual version of Corollary \ref{mininj}, the representation $M$ has a minimal projective resolution of the form
$$0\rightarrow\bigoplus_{i\in I}P_i\otimes{\rm Ext}^1_Q(M,E_i)\rightarrow\bigoplus_{i\in I}P_i\otimes{\rm Hom}_Q(M,E_i)\rightarrow M\rightarrow 0.$$
For all $i\in I$, let $Y_i$ be a vector space of dimension $\langle d,i\rangle-\dim{\rm Hom}_Q(M,E_i)$ (note that this number is non-negative by assumption). Adding the projective representation $\bigoplus_{i\in I}P_i\otimes Y_i$ to the the left and middle term of the above minimal resolution, we claim that the resulting exact sequence $\ses{P'}{P}{M}$ has all the desired properties.\\[1ex]
For any two vertices $i,j\in I$, we have $\langle\dimv P_i,j\rangle=\dim{\rm Hom}_Q(P_i,E_j)=\delta_{i,j}$, and thus $\langle \dimv P,i\rangle=\langle d,i\rangle$ by definition, for all $i\in I$. Since the quiver $Q$ has no oriented cycles, the form $\langle\_,\-\rangle$ is non-degenerate, proving $\dimv P=d$.\\[1ex]
By the assumptions, we have ${\cal M}_{d,n}\not=\emptyset$. Thus, $R_{d,n}^s\not=\emptyset$ is dense in $R_{d,n}$. But the projective representation $P$ has an open orbit in $R_d$, so that its inverse image under the projection $R_{d,n}\rightarrow R_d$ meets $R_{d,n}^s$. Thus, there exists a point $(P,f)\in{\cal M}_{d,n}$. The same argument proves that there exists a point $(P',f')\in {\cal M}_n$.\hb

The above analysis now implies the main theorem of this section.

\begin{theorem}\label{corsurj} The map $\eta:\dot{H}\rightarrow {\bf C}^{A_n}({\cal X}_n)$ is surjective. Its kernel is the linear span of the elements $1_d\chi_{[M]}$ for which there exists an $i\in I$ such that $\langle i,d\rangle>n_i$ or $\dim{\rm Hom}_Q(M,E_i)>\langle d,i\rangle$.
\end{theorem}

\proof The image of the basis element $1_d\chi_{[M]}$ is either zero or a multiple of the characteristic function of the orbit ${\cal O}_{d,[M]}$. The previous proposition now implies the theorem.\hb

\begin{corollary} If $Q$ is of Dynkin type, the algebra ${\bf C}^{A_n}({\cal X}_n)$ is a factor of the modified quantum Borel $\dot{\cal U}_v(\mathfrak{b}^+)$.
\end{corollary}

\proof Using the above theorem and Corollary \ref{raute}, we have in the case of a quiver of Dynkin type a chain of surjective algebra morphisms
$$\dot{\cal U}_v(\mathfrak{b})\stackrel{\sim}{\rightarrow}\dot{C}=\dot{H}\rightarrow{\bf C}^{A_n}({\cal X}_n).$$\hb

\remarks
\begin{enumerate}
\item One can see that, even for fairly simple quivers $Q$, like e.g. $i\rightarrow j\leftarrow k$, there is no obvious extension of the map $\eta$ to a map from the whole modified quantum group $\dot{\cal U}(\mathfrak{g})$ to the convolution algebra ${\bf C}^{A_n}({\cal M}_n\times{\cal M}_n)$.
\item In \cite{BLM}, the convolution algebra on pairs of flags gives a geometric realization of $q$-Schur algebras. One can therefore view the algebras ${\bf C}^{A_n}({\cal X}_n)$ as Borel parts of $q$-Schur algebras ``of Kac-Moody type".
\end{enumerate}

\end{document}